\documentclass[12pt]{article}
\usepackage{amssymb}

\usepackage{makeidx}
\usepackage{amsbsy,amsmath}
\usepackage{latexsym,amssymb,mathptm}
\usepackage{bm}

\newtheorem{Theorem}{Theorem}
\newtheorem{Corollary}[Theorem]{Corollary}
\newtheorem{Lemma}[Theorem]{Lemma}

\begin{document}

\title{On Asymptotic Proximity of Distributions}
\author{{Youri \textsc{Davydov${}^{1}$}}\quad and\quad {\ Vladimir \textsc{Rotar $%
{}^{2}$}}}
\date{{\footnotesize ${}^{1}$ Laboratoire Paul Painlev\'e - UMR 8524}\\
{\footnotesize Universit\'{e} de Lille I - Bat. M2}\\
{\footnotesize \ 59655 Villeneuve d'Ascq, France }\\
{\footnotesize \ Email: youri.davydov@univ-lille1.fr}\\
{\footnotesize \ \vspace{3pt} ${}^{2}$}{\footnotesize \ \ Department of
Mathematics and Statistics }\\
{\footnotesize \ of the San Diego State University, USA and }\\
{\footnotesize \ the Central Economics and Mathematics Institute }\\
{\footnotesize \ of the Russian Academy of Sciences, RF }\\
{\footnotesize \ Email: vrotar@math.ucsd.edu }}
\maketitle


\begin{quote}
\textbf{Abstract.} We consider some general facts concerning convergence
\[
P_{n}-Q_{n}\rightarrow 0\text{\thinspace \thinspace \thinspace as \thinspace
\thinspace \thinspace }n\rightarrow \infty ,
\]
where $P_{n}$ and $Q_{n}$ are probability measures in a complete separable
metric space. The main point is that the sequences $\{P_{n}\}$ and $%
\{Q_{n}\} $ are not assumed to be tight. We compare different possible
definitions of the above convergence, and establish some general properties.



\end{quote}

\section{Introduction and results}

\renewcommand{\theequation}{\thesubsection.\arabic{equation}}

\subsection{Background and Motivation\label{sec.background}}

\setcounter{equation}{0}

Usually, a limit theorem of Probability Theory is a theorem that concerns
convergence of a sequence of distributions $P_{n}$ to a distribution $P$.
However, there is a number of works where the traditional setup is modified,
and the object of study is two sequences of distributions, $\{P_{n}\}\,$and $%
\{Q_{n}\}$. The goal in this case consists in establishing conditions for
convergence
\begin{equation}
P_{n}-Q_{n}\rightarrow 0  \label{i1}
\end{equation}
in a proper sense. In particular problems, $P_{n}$ and $Q_{n}$ are, as a
rule, the distributions of the r.v.'s $f_{n}(X_{1},...,X_{n})$ and $%
f_{n}(Y_{1},...,Y_{n})$, where $f_{n}(\cdot )$ is a function, and $%
X_{1},X_{2},...$ and $Y_{1},Y_{2},...$ are two sequences of r.v.'s. The aim
here is rather to show that different random arguments $X_{1},...,X_{n}$ may
generate close distributions of $f_{n}(X_{1},...,X_{n})$, than to prove that
the distribution of $f_{n}(X_{1},...,X_{n})$ is close to some fixed
distribution (which, above else, may be not true).

Consider, for example, a quadratic form $<A_{n}\mathbf{X}_{n},\mathbf{X}%
_{n}> $, where $A_{n}$ is a $n\times n$-matrix and $\mathbf{X}%
_{n}=(X_{1},...,X_{n})$ is a vector with i.i.d. coordinates. In this case,
the limiting distribution, if any, depends on the matrices $A_{n}$. For
instance, with a proper normalization, the form $%
X_{1}X_{2}+X_{2}X_{3}+...+X_{n-1}X_{n}$ is asymptotically normal (if $%
E\{X_{i}^{2}\}<\infty $), while the form $(X_{1}+...+X_{n})^{2}$ has the $%
\chi _{1}^{2}$ limiting distribution. Nevertheless, one can state the
following unified limit theorem.

Denote by $P_{nF}$ the distribution of $<A_{n}\mathbf{X}_{n},\mathbf{X}_{n}>$
in the case when each $X_{i}$ has a distribution $F$. Then, under rather
mild conditions,
\begin{equation}
P_{nF}-P_{nG}\rightarrow 0  \label{i1+}
\end{equation}
for any two distributions $F$ and $G$ with the same first two moments (see,
\cite{rotar1}, \cite{gotze1} for detail, and references therein). A class $%
\mathcal{F}$ such that (\ref{i1+}) is true for any $F,G\in
\mathcal{F}$ is called an \textit{invariance class} (\cite{rotar1}).

Let us come back to (\ref{i1}). Clearly, such a framework is more general
than the traditional one. First, as was mentioned, the distributions $P_{n}$
and $Q_{n}$ themselves do not have to converge. Secondly, the sequences $%
\{P_{n}\}$ and $\{Q_{n}\}$ are not assumed to be tight, and the convergence
in (\ref{i1}) covers situations when a part of the probability distributions
or the whole distributions ``move away to infinity'' while the distributions
$P_{n}$ and $Q_{n}$ are approaching each other.

To our knowledge, the scheme above was first systematically used in
Lo\'{e}ve \cite[Chapter VIII, Section 28, The Central Asymptotic
Problem]{loeve}, who considered sums of dependent r.v.'s. The same
approach is applied in some non-classical limit theorems for sums of
r.v.'s; that is, theorems not involving the condition of asymptotic
negligibility of separate terms (see, e.g., monograph
\cite{zolotarev} by Zolotarev, survey \cite {rotar1} by Rotar, and
references in \cite{zolotarev} and \cite{rotar1}; Liptser and
Shiryaev \cite{liptser}, and Davydov and Rotar \cite {davydovrotar}
on a non-classical invariance principle.)

Non-linear functions $f_{n}$ have been also considered in the above
framework. In particular, it concerns polynomials, polylinear forms
of r.v.'s, and quasi-polynomial functions; see, e.g., Rotar
\cite{rotar1} and \cite{rotar2} for limit theorems, G\"{o}tze and
Tikhomirov \cite{gotze1}, and Mossel, O'Donnel and Oleszkiewicz
\cite{mossel} for the accuracy of the corresponding invariance
principle. Other interesting schemes different from those above were
explored in D'Aristotile, Diaconis, and Freedman \cite {daristotile}
and in Chatterjee \cite{chatterjee1}.

The present paper addresses general facts on convergence (\ref{i1}).
A corresponding theory was built in Dudley \cite{dudley0},
\cite{dudley}, \cite[Chapter 11]{dudley1} and D'Aristotile,
Diaconis, and Freedman \cite {daristotile}. The paper
\cite{daristotile} concerns some possible definitions of convergence
(\ref{i1}) in terms of uniformities, and establishes connections
between these definitions (see also below). The theory in \cite
{dudley0}-\cite{dudley1} (which is used in part in
\cite{daristotile} also) is mainly based on a metric approach. We
complement and, to a certain extent, develop the theory from
\cite{daristotile} and \cite{dudley1}, paying more attention to a
functional approach. Throughout this paper, we repeatedly refer to
and cite results from \cite{daristotile}-\cite{dudley1}.

First, consider three definitions of convergence (\ref{i1}) explored in \cite
{daristotile} (or, in the terminology of \cite{daristotile}, ``merging'').

\begin{description}
\item[D1.]  $\pi (P_{n},Q_{n})\rightarrow 0$ where $\pi $ is the
L\'{e}vy-Prokhorov metric.

\item[D2.]
\begin{equation}
\int f(x)\left( P_{n}(dx)-Q_{n}(dx)\right) \rightarrow 0  \label{i2}
\end{equation}
for all bounded continuous functions $f$.

\item[D3.]  $T(P_{n})-T(Q_{n})\rightarrow 0$ for all bounded and continuous
(with respect to weak $^{*}$ topology) functions $T$ on the space of
probability measures.
\end{description}

In \cite{daristotile}, it was shown that \textbf{D3} $\Rightarrow$ \textbf{D2} $\Rightarrow$ \textbf{D1}%
, and \textbf{D1,\,D2,\,D3 }are equivalent iff the space on which
measures are defined, is compact. As one can derive from
\cite{daristotile}, and as follows from results below, once we
consider particular sequences $\{P_{n}\}$ and $\{Q_{n}\}$, the above
definitions are equivalent if one of the sequences is tight.

In the general setup, when the above sequences $\{P_{n}\}$ and
$\{Q_{n}\}$ are not assumed to be tight, Definition \textbf{D2} (and
hence \textbf{D3} also) looks too strong. As an example, consider
that from \cite{daristotile}. Let us deal with
distributions on the real line, and let $P_{n}$ be concentrated at the point $n$%
, and $Q_{n}$ -- at the point $n+\frac{1}{n}$. Clearly, (\ref{i2}) is not
true for, say, $f(x)=\sin (x^{2})$. On the other hand, it would have been
unnatural, if a definition had not covered such a trivial case of asymptotic
proximity of distributions. (Clearly, in this case, $\pi
(P_{n},Q_{n})\rightarrow 0$. To make the example simpler, one may consider $%
f(x)=\sin \left( \frac{\pi }{4}x^{2}\right) $. Certainly, the example above
concerns the Euclidean metric in $\Bbb{R}$. For other metrics, points $n$
and $n+\frac{1}{n}$ may be not close. Below, we cover the general case of a
complete separable metric space.)

To fix the situation, one can consider (\ref{i2}) for functions only from
the class of all bounded continuous functions vanishing at infinity, but
such an approach would be too restrictive. In this case, the definition
would not cover situations where parts of the distributions move away to
infinity, continuing to approach each other (or, in the terminology from
\cite{daristotile}, ``merge''). On the other hand, in accordance with the
same definition, the distributions $P_{n}$ and $Q_{n}$ concentrated, for
example, at points $n$ and $2n$, would be viewed as merging, which also does
not look reasonable.

In this paper, we suggest to define weak convergence as that with
respect to all bounded \textit{uniformly }continuous functions.
(This type of convergence was not considered in \cite{daristotile}
and \cite{dudley1}.) We justify this definition proving that such a
convergence is equivalent to convergence in the L\'{e}vy-Prokhorov
metric that satisfies some natural, in our opinion, requirements. In
the counterpart of Definition \textbf{D3}, we require the uniform
continuity of $T$.

We establish also some facts concerning weak convergence uniform on certain
classes of functions $f$ (or linear functionals on the space of
distributions); see Section \ref{sec.uniform} for detail. Proofs turn out to
be, though not very difficult, but not absolutely trivial since the absence
of tightness requires additional constructions. The point is that in the
generalized setup, we cannot choose just one compact, not depending on $n$,
on which all measures will be ``almost concentrated''.

On the other hand, as will follow from results below, if one of the
sequences, $P_{n}$ or $Q_{n}$, is tight, the definition of weak convergence
suggested is equivalent to the classical definition, and we deal with the
classical framework.

We would like to thank P.J.Fitzsimmons and F.D.Lesley for useful discussions.

\subsection{Main Results}

\setcounter{equation}{0}

\subsubsection{Weak convergence and the L\'{e}vy-Prokhorov metric\label%
{sec.weak}}

Let $(\Bbb{H},\,\rho )$ be a complete separable metric space, and $\mathcal{B%
}$ be the corresponding Borel $\sigma $-algebra.

The symbols \ $P$ and $Q$, with or without indices, will denote probability
distributions on $\mathcal{B}$. All functions $f$ below, perhaps with
indices, are continuous functions $f:\Bbb{H\rightarrow R}$.

We denote by $\mathcal{C}$ the class of all bounded and continuous functions
on $(\Bbb{H},\,\rho )$, and by $\overline{\mathcal{C}}$ - the class of all
bounded and uniformly continuous functions.

For two sequences of probability measures (distributions), $\{P_{n}\}$ and $%
\{Q_{n}\}$, we say that $P_{n}-Q_{n}\rightarrow 0\,\,\,\,\,\,\,$\textit{%
weakly with respect to (w.r.t.) a class of functions }$\mathcal{K}$ if $\int
fd(P_{n}-Q_{n})\rightarrow 0$ for all\textit{\ }$f\in \mathcal{K}$\textit{. }

If we do not mention a particular class $\mathcal{K}$, \textit{the term weak
convergence (or more precisely, merging) will concern that w.r.t.} $%
\overline{\mathcal{C}}$. When it cannot cause misunderstanding,  we will use
the term ``convergence'' in the situation of merging also.

In the space of probability distributions on $\mathcal{B}$, we define - in a
usual way - the \textit{L\'{e}vy-Prokhorov metric}
\begin{equation}
\pi (P,Q)=\inf \{\varepsilon :P(A^{\varepsilon })\leq Q(A)+\varepsilon \text{
\thinspace \thinspace for all closed sets }A\}.  \label{i3}
\end{equation}
(One can restrict himself to only one inequality, not switching $P$ and $Q$;
see Dudley \cite[Theorem 11.3.1]{dudley}.)

Our main result is

\begin{Theorem}
\label{th1}The difference
\begin{equation}
P_{n}-Q_{n}\rightarrow 0\,\,\,\,\,\,\,\,\text{weakly (w.r.t. }\overline{%
\mathcal{C}}\text{) }  \label{new4}
\end{equation}
\textit{\ if and only if }
\begin{equation}
\pi (P_{n},Q_{n})\rightarrow 0.  \label{new3}
\end{equation}
\end{Theorem}

The choice of the L\'{e}vy-Prokhorov metric as a ``good'' metric that
justifies the definition of weak convergence above, is connected, first of
all, with the fact that the analog of the Skorokhod representation theorem (%
\cite{skorohod}, see also, e.g., \cite[Sec.11.7]{dudley1}) holds in the case
of merging measures. More precisely, the following is true.

Let, the symbols $X$ and $Y$, perhaps with indices, denote random variables
defined on a probability space $(\Omega ,\mathcal{A},P)$, and assuming
values in $\Bbb{H}$ (that is, these variables are $\mathcal{A}\rightarrow
\mathcal{B}$ measurable.) The symbol $P_{X}$ stands for the distribution of $%
X$.

Below, two metrics $r_{1}(x,y)$ and $r_{2}(x,y)$ in a space are said
to
be \textit{uniformly equivalent}, if for any two sequences $\{x_{n}\}$ and $\{y_{n}\}$%
, the relations $r_{1}(x_{n},y_{n})\rightarrow 0$ and $r_{2}(x_{n},y_{n})%
\rightarrow 0$ are equivalent.

The first two (and main) assertions of the next theorem are stated and
proved in Dudley \cite{dudley}; see also the second edition
\cite[Theorem 11.7.1]{dudley1}. For the completeness of the picture, we
present all facts regarding the metric $\pi$ in one theorem.

\begin{Theorem}
\label{th2} Metric $\pi $ is the only metric, to within uniform equivalence,
that satisfies the following conditions.

\begin{enumerate}
\item[A.]  If $\rho (X_{n},Y_{n})\stackrel{P}{\rightarrow }0$, then $\pi
(P_{X_{n}},\,P_{Y_{n}})\rightarrow 0$.

\item[B.]  If $\pi (P_{n},\,Q_{n})\rightarrow 0$, then there exist a
probability space and random elements $X_{n},Y_{n}$ on this space such that $%
P_{X_{n}}=P_{n}$, $Q_{n}=P_{Y_{n}}$, and $\rho (X_{n},Y_{n})\stackrel{P}{%
\rightarrow }0$.

\item[C.]  If $\pi (P_{n},\,Q_{n})\rightarrow 0$, then $\pi (P_{n}\circ
f^{-1},Q_{n}\circ f^{-1})\rightarrow 0$ for any uniformly continuous
function $f$.

\item[D.]  If $Q_{n}=Q$, then the convergence $\pi (P_{n},\,Q)\rightarrow 0$
is equivalent to weak convergence with respect to all bounded continuous
functions.$\vspace{0.1in}$
\end{enumerate}
\end{Theorem}

\textbf{Remarks.}

\begin{enumerate}
\item  We have already mentioned above an example showing that convergence (%
\ref{i2}) for all $f\in \mathcal{C}$ does not possess Property A. In other
words, there exist r.v.'s $X_{n}$ and $Y_{n}$ such that $\rho (X_{n},Y_{n})%
\stackrel{P}{\rightarrow }0$, while $P_{X_{n}}-P_{Y_{n}}\nrightarrow 0\,\,\,$%
weakly with respect to $\mathcal{C}$. If $\Bbb{H}$ has bounded but
non-compact sets, a similar example may be constructed even for continuous
functions equal to zero out of a closed bounded set.

Indeed, let $O_{r}(x)=\{y:\rho (x,y)\leq r\}$. Consider the case where for
some $x_{0}$, the ball $O=O_{1}(x_{0})$ is not compact. Then there exists a
sequence $\{x_{n}\}\subset O$ which does not contain a converging
subsequence. Furthermore, there exists a numerical sequence $\delta
_{n}\rightarrow 0$, such that the balls $O_{\delta _{n}}(x_{n})$ are
disjoint, and each contains only one element from $\{x_{n}\}$, that is, the
center $x_{n}$. We define
\[
f(x)=\left\{
\begin{array}{ll}
1-\frac{1}{\delta _{n}}\rho (x,x_{n}) & \text{if }x\in O_{\delta
_{n}}(x_{n}), \\
0 & \text{if }x\not{\in}\cup _{k}O_{\delta _{k}}(x_{k}).
\end{array}
\right.
\]
The function $f$ is, first, bounded, and secondly, due to the choice
of $\{x_{n}\}$, $f$ is continuous. On the other hand, one may set $%
X_{n}\equiv x_{n} $ and $Y_{n}\equiv y_{n}$, where $y_{n}$ is a point from
the boundary of $O_{\delta _{n}}(x_{n})$. Clearly, $\rho
(X_{n},Y_{n})=\delta _{n}\rightarrow 0$, while $\int fdP_{X_{n}}=1$ and $%
\int fdP_{Y_{n}}=0$.

\item  As was mentioned in the introduction, if one of the sequences, say $%
\{Q_{n}\}$, is tight, and (\ref{new4}) is true, then the other sequence, $%
\{P_{n}\}$, is also tight. In this case, relation (\ref{i2}) is true for all
$f\in \mathcal{C}$, and we deal with the classical scheme. If $(\Bbb{H},\rho
)$ is a space in which any closed bounded set is compact, this fact is easy
to prove directly. In the general case, it is easier to appeal to Theorem
\ref{th1} and Prokhorov's theorem on relative compactness w.r.t. functions
from $\mathcal{C}$ (\cite{prokhorov}).

More precisely, assume that (\ref{new4}) holds and $\{Q_{n}\}\,$is tight.
Then, by Pro\-kho\-rov's theorem, $\{Q_{n}\}$ is relatively compact with
respect to weak convergence for functions from $\mathcal{C}$. Let a
subsequence $Q_{n_{k}}$ weakly converges to some $Q$ w.r.t. $\mathcal{C}$,
and hence $\pi (Q_{n_{k}},Q)\rightarrow 0$. By virtue of Theorem \ref{th1}, $%
\pi (P_{n_{k}},Q_{n_{k}})\rightarrow 0$, and consequently, $\pi
(P_{n_{k}},Q)\rightarrow 0$. So, $P_{n_{k}}$ converges to the same $Q$
w.r.t. to $\mathcal{C}$, and
\begin{equation}
P_{n_{k}}-Q_{n_{k}}\rightarrow 0\text{ w.r.t. to all functions from }%
\mathcal{C}\text{.}  \label{new5}
\end{equation}
Thus, any subsequence of $\{P_{n}\}$ contains a subsubsequence convergent
w.r.t. $\mathcal{C}$. Hence, again by Prokhorov's theorem, $\{P_{n}\}$ is
tight.

Moreover, in this case $P_{n}-Q_{n}\rightarrow 0$ w.r.t. $\mathcal{C}$.
Indeed, otherwise we could select convergent subsequences $P_{n_{k}}$ and $%
Q_{n_{k}}$ and a bounded and continuous $f$ such that $\int
fd(P_{n_{k}}-Q_{n_{k}})\nrightarrow 0$, which would have contradicted (\ref
{new5}).

\item  Let us return to Definition \textbf{D3} from Section \ref{sec.background}. To
make it suitable to the setup of this paper, one should consider functions $T
$ on the space of probability measures on $\mathcal{B}$, \textit{uniformly}
continuous with respect to the the L\'{e}vy-Prokhorov metric (or, which is
the same, w.r.t. the weak convergence regarding $\overline{\mathcal{C}}$).
So modified Definition \textbf{D3} will be equivalent to convergences (\ref{new4})-(%
\ref{new3}).
\end{enumerate}

\subsubsection{Lipschitz functions}

For $f\in \mathcal{C}$, we set
\[
\left\| f\right\| _{L}=\sup_{x\neq y}\left\{ \frac{|f(x)-f(y)|}{\rho (x,y)}%
\right\} ,
\]
and $\mathcal{C}_{BL}=\{f:\left\| f\right\| _{L}<\infty ,\,\left\| f\right\|
_{\infty }<\infty \}$.

In Dudley \cite[Sec.11.3]{dudley1}, it was shown that in the traditional
setup, for the weak convergence $P_{n}\rightarrow P$ w.r.t. the whole class $%
\mathcal{C}$, it suffices to have the convergence w.r.t. $\mathcal{C}_{BL}$.
A similar property is true for the generalized setup of merging probability
measures. Namely, let $\{P_{n}\}$ and $\{Q_{n}\}$ be two fixed sequences of
probability measures.

\begin{Theorem}
\label{th3}\thinspace \thinspace \thinspace \thinspace \thinspace The weak
convergence $P_{n}-Q_{n}\rightarrow 0$ with respect to $\mathcal{C}_{BL}$%
\thinspace implies the weak convergence with respect to $\overline{\mathcal{C%
}}$.
\end{Theorem}

Next, we consider classes of functions on which weak convergence is uniform.

\subsubsection{Uniform convergence\label{sec.uniform}}

In Dudley \cite{dudley1}, it was shown that in the traditional setup, weak
convergence is uniform on any class of functions $f$ with uniformly bounded
norms $\left\| f\right\| _{L}$ and $\,\left\| f\right\| _{\infty }$. More
precisely, consider the metric
\[
\beta (P,Q):=\sup \left\{ \left| \int fd(P-Q)\right|
:\,\,\,||f||_{L}+||f||_{\infty }\leq 1\right\}.
\]
In \cite[Section 1.3]{dudley1}, it was proved that the weak convergence $%
P_{n}\rightarrow P$ w. r. t. $\mathcal{C}$, is equivalent to the convergence
$\beta (P_{n},P)\rightarrow 0$.

We establish a similar property in the generalized setup and for
arbitrary classes of functions with a fixed order of their moduli of
continuity.

For $f\in \overline{\mathcal{C}}$, we define its modulus of continuity
\begin{equation}
\omega _{f}(h)=\sup \{|f(x)-f(y)|:\rho (x,y)\leq h\}.  \label{new01}
\end{equation}

Let $\omega (h)$ be a \textit{fixed} non-decreasing function on $\Bbb{R}^{+}$%
, such that $\omega (h)\rightarrow 0$ as $h\rightarrow 0$. Set
\[
\mathcal{C}_{\omega }=\{f:\left\| f\right\| _{\infty }<\infty ,\text{
\thinspace \thinspace and \thinspace \thinspace }\omega _{f}(h)=O(\omega
(h)+h)\}.
\]
(Usually, $h=O(\omega(h))$, and one can write just $\omega _{f}(h)=O(\omega
(h))$. However, there are situations when $\omega(h)$ even equals zero for
sufficiently small $h$'s; for example, if $\Bbb{H}=\Bbb{N}$ with the usual
metric.)

The next proposition, having its intrinsic value, plays an essential role in
proving Theorem \ref{th1}.

\begin{Theorem}
\label{th4}Let
\begin{equation}
\int fd\left( P_{n}-Q_{n}\right) \rightarrow 0  \label{new02}
\end{equation}
for all $f\in \mathcal{C}_{\omega }$, and let
\[
\mathcal{F}_{\omega }=\{f:\left\| f\right\| _{\infty }<1\text{, \thinspace
\thinspace and\thinspace \thinspace \thinspace \thinspace }\omega
_{f}(h)\leq \omega (h)\text{ \thinspace \thinspace for all \thinspace
\thinspace }h\geq 0\}.
\]
Then
\begin{equation}
\sup_{f\in \mathcal{F}_{\omega }}\left| \int fd\left( P_{n}-Q_{n}\right)
\right| \rightarrow 0.  \label{new03}
\end{equation}
\end{Theorem}

Clearly, instead of the above class $\mathcal{F}_{\omega }$, one may
consider a class of uniformly bounded (not necessarily by one)
functions $f$ such that $\omega _{f}(h)\leq k_{1}\omega (h)+k_{2}h$,
where $k_{1},\,k_{2}$
are fixed constants. (Such a formal generalization follows from (\ref{new03}%
) just by replacing $\omega (h)$ by $k_{1}\omega (h)+k_{2}h$.)

\begin{Corollary}
\label{cor1}If (\ref{new02}) is true for all $f\in \mathcal{C}_{BL}$, then
\begin{equation}
\sup_{f\in \mathcal{F}}\left| \int fd\left( P_{n}-Q_{n}\right) \right|
\rightarrow 0  \label{new04}
\end{equation}
for any class $\mathcal{F}$ of uniformly bounded functions with uniformly
bounded Lipschitz constants.
\end{Corollary}

\begin{Corollary}
\label{cor2}If (\ref{new02}) is true for all $f\in \overline{\mathcal{C}}$,
then (\ref{new04}) is true for any class $\mathcal{F}$ of uniformly bounded
and uniformly equicontinuous functions.
\end{Corollary}

To derive the above corollaries from Theorem \ref{th4}, one should set $%
\omega (h)=\newline \omega _{\mathcal{F}}(h):=\sup_{f\in
\mathcal{F}}\omega _{f}(h)$. (In the literature, there is no unity
in definitions of equicontinuity: some authors define it pointwise;
in other definitions, the word ``uniformly'' is redundant. When
talking about uniform equicontinuity, we mean that the function
$\omega _{\mathcal{F}}(h)$ is bounded and vanishing at the origin.)

\section{Proofs}

\setcounter{equation}{0}

Proofs of Theorems \ref{th1} and \ref{th3} essentially use Corollary \ref
{cor1} from Theorem \ref{th4} and Theorem \ref{th2}. We start with a proof
of the latter theorem -- the main assertions of this theorem, A and B, are
known (\cite[Sec.11.7]{dudley1}), and the rest of the proof is short. The
proof of Theorem \ref{th4} is relegated to the last Section \ref{sec.sup}.

\subsection{Proofs of Theorems \ref{th1} -- \ref{th3}}

\setcounter{equation}{0}

\subsubsection{Proof of Theorem \ref{th2}}

To justify Property C, consider the r.v.'s $X_{n},Y_{n}$ defined in Property
B. We have $|f(X_{n})-f(Y_{n})|\leq \omega _{f}(\rho (X_{n},Y_{n}))$.
Consequently, $f(X_{n})-f(Y_{n})\stackrel{P}{\rightarrow }0$, and it
suffices to appeal to Property A.

Property D is obvious. Now, let $r_{1}(\cdot ,\cdot )$ and $r_{2}(\cdot
,\cdot )$ \ be two metrics with Properties A-B. If $r_{1}(P_{n},Q\,_{n})%
\rightarrow 0$, then there exist $X_{n},\,Y_{n}$ for which $\rho
(X_{n},Y_{n})\rightarrow 0$, and hence by Property A, $r_{2}(P_{n},Q\,_{n})%
\rightarrow 0$. \thinspace \thinspace $\blacksquare \vspace{0.1in}$

For proving Theorems \ref{th1} and \ref{th3}, we need

\subsubsection{Two lemmas}

\begin{Lemma}
\label{lemma2}If $\pi (P_{n},Q\,_{n})\rightarrow 0$, then
\begin{equation}
\sup_{f\in \mathcal{F}}\int fd(P_{n}-Q_{n})\rightarrow 0  \label{new05}
\end{equation}
for any class $\mathcal{F}$ of uniformly bounded and uniformly
equicontinuous functions.
\end{Lemma}

\textbf{Proof. }By Theorem \ref{th2}, there exist $X_{n},Y_{n}$ such that $%
P_{X_{n}}=P_{n}$, $Q_{n}=P_{Y_{n}}$, and $\rho (X_{n},Y_{n})\stackrel{P}{%
\rightarrow }0$. \ By conditions of the lemma, $M:=\sup_{f\in \mathcal{F}%
}\left\| f\right\| _{\infty }<\infty $, and $\omega (h):=\sup_{f\in \mathcal{%
F}}\omega _{f}(h)\rightarrow 0$ as $h\rightarrow 0$. For any $\varepsilon >0$%
,
\begin{eqnarray*}
\left| \int fd(P_{n}-Q_{n})\right| &=&\left| E\left\{
f(X_{n})-f(Y_{n})\right\} \right| \leq E\left\{ \left|
f(X_{n})-f(Y_{n})\right| \right\} \\
&=&E\{\left| f(X_{n})-f(Y_{n})\right| ;\,\rho (X_{n},Y_{n})>\varepsilon \} \\
&+&E\{\left| f(X_{n})-f(Y_{n})\right| ;\,\rho (X_{n},Y_{n})\leq \varepsilon
\} \\
&\leq &2MP(\rho (X_{n},Y_{n})>\varepsilon )+\omega (\varepsilon ).
\end{eqnarray*}
Hence,
\[
\overline{\lim_{n}}\sup_{f\in \mathcal{F}}\left| \int fd(P_{n}-Q_{n})\right|
\leq \omega (\varepsilon )\rightarrow 0\text{ as }\varepsilon \rightarrow
0.\,\,\,\,\,\,\,\,\,\,\,\,\,\blacksquare
\]

\begin{Lemma}
\label{lemma2+}(Dudley \cite[a part of Theorem 11.7.1]{dudley1}). Suppose (%
\ref{new05}) is true for $\mathcal{F}=\mathcal{F}_{1}:=\{f:\left\| f\right\|
_{\infty }\leq 1,\,\,\left\| f\right\| _{L}\leq 1\}$. Then $\pi
(P_{n},Q\,_{n})\rightarrow 0$.
\end{Lemma}

\textbf{Proof.} In \cite{dudley1}, the proof of this fact is based on the
relation $\pi \leq 2\sqrt{\beta }$, which is proved separately. For the
completeness of the picture, we give a direct proof (which, in essence, is
very close to the reasoning in \cite[p.396 ]{dudley1}).

For a closed set $K$, and an $\varepsilon >0$, set
\[
I_{K}^{\varepsilon }(x)=\left\{
\begin{array}{ll}
1 & \text{if }x\in K, \\
1-\frac{1}{\varepsilon }\rho (x,K) & \text{if }x\in K^{\varepsilon
}\setminus K, \\
0 & \text{otherwise.}
\end{array}
\right.
\]
(Here, $K^{\varepsilon }$ is the $\varepsilon $-neighborhood of $K$.)

Since $\omega _{I_{K}^{\varepsilon }}(h)\leq \frac{1}{\varepsilon }h$, the
family $\{I_{K}^{\varepsilon }(x):K$ is closed$\}\subset \mathcal{F}%
_{\varepsilon }:=\{f:\left\| f\right\| _{\infty }\leq 1,\,\,\left\|
f\right\| _{L}\leq 1/\varepsilon \}$. Clearly, if (\ref{new05}) holds for $%
\mathcal{F}=\mathcal{F}_{1}$, then it holds for $\mathcal{F}=\mathcal{F}%
_{\varepsilon }$ for any $\varepsilon >0$. Therefore,
\begin{equation}
\Delta _{n}(\varepsilon ):=\sup_{K}\left| \int I_{K}^{\varepsilon
}(x)d(P_{n}-Q_{n})\right| \rightarrow 0  \label{new06}
\end{equation}
for any $\varepsilon >0$ as $n\rightarrow \infty $.

On the other hand,
\begin{equation}
P_{n}(K)\leq \int I_{K}^{\varepsilon }(x)dP_{n}\leq \Delta _{n}(\varepsilon
)+Q_{n}(K^{\varepsilon }).  \label{new2}
\end{equation}
From (\ref{new06}) and (\ref{new2}), it follows that for sufficiently large $%
n$,
\[
P_{n}(K)\leq \varepsilon +Q_{n}(K^{\varepsilon }),
\]
which implies that $\pi (P_{n},Q_{n})\leq \varepsilon $. Since $\varepsilon $
is arbitrary small, this means that $\pi (P_{n},Q_{n})\rightarrow 0$. $%
\blacksquare $

\subsubsection{Proof of Theorem \ref{th1}\label{sec.proof1}}

The implication $\pi (P_{n},Q_{n})\rightarrow 0$ $\Rightarrow $ $\int
fd(P_{n}-Q_{n})\rightarrow 0$ for any $f\in \overline{\mathcal{C}}$,
immediately follows from Lemma \ref{lemma2}.

Assume $\int fd(P_{n}-Q_{n})\rightarrow 0$ for any $f\in \overline{\mathcal{C%
}}$. Then, the same is true for all $f\in \mathcal{C}_{BL}$.

Hence, by Corollary \ref{cor1} from Theorem \ref{th4}, relation (\ref{new05}%
) holds for $\mathcal{F}=\mathcal{F}_{1}$ (defined in Lemma
\ref{lemma2+}). This implies the convergence $\pi
(P_{n},Q\,_{n})\rightarrow 0$ by virtue of Lemma
\ref{lemma2+}.\thinspace \thinspace \thinspace $\blacksquare $

\subsubsection{Proof of Theorem \ref{th3}}

Let $\int fd(P_{n}-Q_{n})\rightarrow 0$ for any $f\in \mathcal{C}_{BL}$.

As was proved in Section \ref{sec.proof1} above, $\pi
(P_{n},Q\,_{n})\rightarrow 0$. By virtue of Theorem \ref{th1}, this implies
that $\int fd(P_{n}-Q_{n})\rightarrow 0$ for all $f\in \overline{\mathcal{C}}
$. \thinspace \thinspace $\blacksquare $

\subsection{Proof of Theorem \ref{th4}\label{sec.sup}}

\setcounter{equation}{0}

\subsubsection{Three more lemmas}

\begin{Lemma}
\label{lemma5}For any two functions $f$ and $g$,
\begin{eqnarray}
\omega _{fg}(h) &\leq &||g||_{\infty }\omega _{f}(h)+||f||_{\infty }\,\omega
_{g}(h),  \label{new071} \\
\omega _{f\vee g}(h) &\leq &\max \{\omega _{f}(h),\,\omega _{g}(h)\},
\label{new072} \\
\omega _{f\wedge g}(h) &\leq &\max \{\omega _{f}(h),\,\omega _{g}(h)\}.
\label{new073}
\end{eqnarray}
provided that the l.-h.sides are finite.
\end{Lemma}

\textbf{Proof }is straightforward and very close to that in
\cite[Propositions 11.2.1-2]{dudley1} dealing with Lipschitz functions. $%
\vspace{0.1in}$

Next, for a function $f(x)$, a set $K$, and a number $t>0$, we define the
function
\begin{equation}
f_{K}^{(t)}(x)=\left\{
\begin{array}{ll}
f(x) & \text{if }x\in K, \\
f(x)\left( 1-\frac{\rho (x,K)}{t}\right) & \text{ if }x\in K^{t}\setminus K,
\\
0 & \text{otherwise.}
\end{array}
\right.  \label{f-t}
\end{equation}

\begin{Lemma}
\label{lemma3}Let $f$ \thinspace \thinspace be a bounded uniformly
continuous function. Then
\begin{equation}
\left\| f_{K}^{(t)}\right\| _{\infty }\leq \left\| f\right\| _{\infty },
\label{l31}
\end{equation}
and for any $h\geq 0$,
\begin{equation}
\omega _{f_{K}^{(t)}}(h)\leq \omega _{f}(h)+\left\| f\right\| _{\infty }%
\frac{h}{t}.  \label{l32}
\end{equation}
\end{Lemma}

\textbf{Proof.} Bound (\ref{l31}) is obvious. Next, note that for any $x\,$%
and $y$,
\begin{equation}
\left| \rho (x,K)-\rho (y,K)\right| \leq \rho (x,y).  \label{l33}
\end{equation}
(Indeed, for any $z\in K$,
\[
\rho (x,K)\leq \rho (x,z)\leq \rho (x,y)+\rho (y,z),
\]
which implies that $\rho (x,K)\leq \rho (x,y)+\rho (y,K)$. We can also
switch $x$ and $y$.)

Thus, for $q_{t}(x):=1-\rho (x,K)/t$, we have
\[
\omega _{q_{t}}(h)\leq h/t.
\]

Together with (\ref{new071}), this implies (\ref{l32}). \thinspace
\thinspace $\blacksquare $\vspace{0.1in}

Note that, in particular, from Lemma \ref{lemma3} it follows that if $f$ is
uniformly continuous, so does $f_{K}^{(t)}$.

Let the signed measure $\Psi _{n}=P_{n}-Q_{n}$.

\begin{Lemma}
\label{lemma4}Suppose
\begin{equation}
\int fd\Psi _{n}\rightarrow 0  \label{p01}
\end{equation}
for any $f\in \overline{\mathcal{C}}$. Let $\mathcal{F}$ be a class of
uniformly bounded and uniformly equicontinuous functions. Set
\begin{equation}
\omega (h):=\sup_{f\in \mathcal{F}}\omega _{f}(h).  \label{p0-}
\end{equation}
Then for any compact $K$, and $t>0$,
\begin{equation}
\overline{\lim_{n}}\sup_{f\in \mathcal{F}}\left| \int f_{K}^{(t)}d\Psi
_{n}\right| \leq 4\omega (t).  \label{p0}
\end{equation}
\end{Lemma}

\textbf{Proof.} By the Arzel\`{a}-Ascoli theorem, for any $\varepsilon >0$,
there is a finite family $\{f_{1},...,f_{d}\}\subset \mathcal{F}$ such that
for any $f\in \mathcal{F}$, there exists $s=s(f)\in \{1,...,d\}$, for which
\begin{equation}
\sup_{x\in K}|f(x)-f_{s}(x)|<\varepsilon .  \label{p1}
\end{equation}
On the other hand, for any $z\in K$, and $x\in K^{t}$,
\begin{eqnarray*}
|f_{K}^{(t)}(x)-f_{s,K}^{(t)}(x)| &=&|f(x)-f_{s}(x)|\left( 1-\frac{\rho (x,K)%
}{t}\right) \\
&\leq &|f(x)-f(z)|+|f_{s}(x)-f_{s}(z)|+|f(z)-f_{s}(z)| \\
&\leq &2\omega (\rho (x,z))+\varepsilon .
\end{eqnarray*}
Hence, for $x\in K^{t}$,
\begin{equation}
|f_{K}^{(t)}(x)-f_{s,K}^{(t)}(x)|\leq 2\omega (t)+\varepsilon .  \label{p2}
\end{equation}

Therefore,
\begin{eqnarray*}
\left| \int f_{K}^{(t)}d\Psi _{n}\right| &\leq &\left| \int
(f_{K}^{(t)}(x)-f_{s,K}^{(t)}(x))d\Psi _{n}\right| +\left| \int
f_{s,K}^{(t)}(x)d\Psi _{n}\right| \\
&\leq &2(2\omega (t)+\varepsilon )+\left| \int f_{s,K}^{(t)}(x)d\Psi
_{n}\right| \\
&\leq &4\omega (t)+2\varepsilon +\max_{m\in \{1,...,d\}}\left| \int
f_{m,K}^{(t)}(x)d\Psi _{n}\right| .
\end{eqnarray*}
Since
\begin{equation}
\max_{m\in \{1,...,d\}}\left| \int f_{m,K}^{(t)}(x)d\Psi _{n}\right|
\rightarrow 0\text{ as }n\rightarrow \infty ,  \label{p3}
\end{equation}
and $\varepsilon $ is arbitrary small, this implies (\ref{p0}). $%
\blacksquare $ \vspace{.1in}

We turn to

\subsubsection{A direct proof of Theorem \ref{th4}}

Consider a fixed function $\omega (h)$, and the class $\mathcal{F}=\mathcal{F%
}_{\omega }$ from the statement of Theorem \ref{th4}. Assume that there
exist a sequence $\{f_{1},f_{2},...\,\,\}\subset \mathcal{F}$, a sequence $%
\{m_{k}\}$, and a $\delta >0$, such that $\left| \int f_{k}d\Psi
_{m_{k}}\right| \geq \delta $\thinspace \thinspace \thinspace for
all $k=1,2,...\,\,$. Without loss of generality we can identify
$\{m_{k}\}$ with $\Bbb{N}$, and write just
\[
\left| \int f_{k}d\Psi _{k}\right| \geq \delta \text{\thinspace
\thinspace \thinspace for all }k=1,2,...\,\,.
\]
We may also assume that for all $k$'s,
\begin{equation}
0\leq f_{k}(x)\leq 1.  \label{p3+}
\end{equation}

Now, let $t$ and $\varepsilon \leq t$ be two \textit{fixed} positive numbers
to be chosen later. Let $n_{1}=1$, and $K_{1}$ be a compact such that
\[
|\Psi _{n_{1}}|(K_{1}^{\complement})\leq t.
\]
(Here the measure $|\Psi |(dx)$ is the variation of $\Psi (dx)$, and $%
K^{\complement}$ is the complement of $K$.)

Let $n_{2}\geq n_{1}+1$ be a number such that for all $n\geq n_{2}$%
\begin{equation}
\left| \int f_{1,K_{1}}^{(\varepsilon /3)}d\Psi _{n}\right| \leq t,
\label{p5}
\end{equation}
and
\begin{equation}
\sup_{j}\left| \int f_{j,K_{1}}^{(t)}d\Psi _{n}\right| \leq 5\omega
(t). \label{p6}
\end{equation}
Inequality (\ref{p5}) is true for sufficiently large $n$ because
$f_{1,K_{1}}^{(\varepsilon /3)}$ is uniformly continuous due to
Lemma \ref{lemma3}, and (\ref{p6}) holds for large $n$ by virtue of
(\ref{p0}). Note also that in (\ref{p6}) we deal with
the supremum over a class of functions, while (\ref{p5}) concerns \textit{%
only one fixed} function.

Next, we consider a compact $K_{2}$ such that $K_{2}\supset K_{1}$, and
\[
|\Psi _{n_{2}}|(K_{2}^{\complement})\leq t.
\]

Let us set $\tilde{f}_{1}=f_{1}$, and
\[
\tilde{f}_{2}=f_{n_{2}}-f_{n_{2},K_{1}}^{(t)}.
\]
By virtue of (\ref{p6}),
\[
\int \tilde{f}_{2}d\Psi _{n_{2}}=\int f_{n_{2}}d\Psi _{n_{2}}-\int
f_{n_{2},K_{1}}^{(t)}d\Psi _{n_{2}}\geq \delta -5\omega (t).
\]
Also,
\[
\Vert \tilde{f}_{2}\Vert _{\infty }\leq 1,\,\,\,\text{and }\,\,\,\,\tilde{f}%
_{2}(x)=0\text{ \thinspace for all }x\in K_{1}\text{.}
\]
By Lemma \ref{lemma3},
\begin{equation}
\omega _{\tilde{f}_{2}}(h)\leq \omega _{f_{n_{2}}}(h)+\omega
_{f_{n_{2},K_{1}}^{(t)}}(h)\leq 2\omega (h)+\frac{h}{t}.  \label{p611}
\end{equation}

Now, we set $L_{1}=K_{1}$, and $L_{2}=K_{2}\setminus K_{1}^{\varepsilon }$.
Let
\begin{equation}
g_{1}(x)=\tilde{f}_{1,K_{1}}^{(\varepsilon /3)}(x)\text{, \thinspace
\thinspace and \thinspace \thinspace \thinspace }g_{2}(x)=g_{1}(x)+\tilde{f}%
_{2,L_{2}}^{(\varepsilon /3)}(x).  \label{p61}
\end{equation}
By construction,
\[
g_{1}(x)=0\text{ \thinspace \thinspace for }x\notin K_{1}^{\varepsilon /3}%
\text{, and \thinspace \thinspace }g_{2}(x)=0\,\,\,\,\,\text{for }x\notin
K_{2}^{\varepsilon /3}.
\]
By Lemma \ref{lemma3},
\begin{equation}
\omega _{g_{1}}(h)\leq \omega (h)+\frac{h}{\varepsilon /3}=\omega (h)+\frac{%
3h}{\varepsilon }.  \label{p163}
\end{equation}
Now, since the sets $L_{1}^{\varepsilon /3}$ and $L_{2}^{\varepsilon /3}$
are disjoint, in (\ref{p61}), either $g_{1}(x)$ or $\tilde{f}%
_{2,L_{2}}^{(\varepsilon /3)}(x)\,\,\,$equals zero. So, we can also write
that $g_{2}(x)=\max \{g_{1}(x),\,\,\tilde{f}_{2,L_{2}}^{(\varepsilon
/3)}(x)\}$. Then, from Lemma \ref{lemma5}, Lemma \ref{lemma3}, (\ref{p163}),
and (\ref{p611}), it follows that
\begin{equation}
\omega _{g_{_{2}}}(h)\leq \max \{\omega _{g_{1}}(x),\,\omega _{\tilde{f}%
_{2}}(h)+\frac{h}{\varepsilon /3}\}\leq 2\omega (h)+\frac{h}{t}+\frac{3h}{%
\varepsilon }.  \label{p612}
\end{equation}
In view of (\ref{p163}), bound (\ref{p612}) is true for both functions, $%
g_{1}$ and $g_{2}$.

Thus, both functions, $g_{1}$ and $g_{2}$, are bounded and uniformly
continuous (since $\varepsilon ,t>0$ are fixed).

Next, we choose $n_{3}\geq n_{2}+1$ such that for all $n\geq n_{3}$,
\[
\left| \int g_{2}(x)d\Psi _{n}\right| \leq t,
\]
and
\[
\sup_{j}\left| \int f_{j,K_{2}}^{(t)}d\Psi _{n}\right| \leq 5\omega (t).
\]
Let $K_{3}$ be a compact such that $K_{3}\supset K_{2}$, and
\[
|\Psi _{n_{3}}|(K_{3}^{\complement })\leq t.
\]
We define a function
\[
\tilde{f}_{3}=f_{n_{3}}-f_{n_{3},K_{2}}^{(t)}
\]
which has properties similar to those of\thinspace \thinspace $\tilde{f}_{2}$%
, and we define the set
\[
L_{3}=K_{3}\setminus K_{2}^{\varepsilon }.
\]
The sets $L_{1}^{\varepsilon /3}$, $L_{2}^{\varepsilon /3}$, and $%
L_{3}^{\varepsilon /3}$ are mutually disjoint. We set
\begin{equation}
g_{3}(x)=g_{2}(x)+\tilde{f}_{3,L_{3}}^{(\varepsilon /3)}(x),  \label{p62}
\end{equation}
and again note that in (\ref{p62}), either $g_{2}(x)$ or $\tilde{f}%
_{2,L_{3}}^{(\varepsilon /3)}(x)$ $\,$equals zero. Similarly to what we did
above, we conclude that (\ref{p612}) is true for $g_{3}$ also. So, $g_{3}(x)$
is fixed, uniformly continuous, and $g_{3}(x)=0$ for $x\not%
{\in}K_{3}^{\varepsilon /3}$.

Continuing the recurrence procedure in the same fashion, we come to the
following objects.

\begin{description}
\item[(a)]  The sequence $n_{m}\rightarrow \infty $.

\item[(b)]  The sequence of compacts $K_{m}$ such that for all $m$%
\begin{equation}
|\Psi _{n_{m}}|(K_{m}^{\complement })\leq t.  \label{n0}
\end{equation}

\item[(c)]  The sequence of compact sets $L_{m}\subset K_{m}$ such that the
sets $L_{m}^{\varepsilon /3}$ are disjoint.

\item[(d)]  The sequence of functions $\tilde{f}_{m}$ such that for all $m$%
\begin{equation}
\tilde{f}_{m}(x)=0\text{ \thinspace for all }x\in K_{m-1}\text{,}
\label{n11}
\end{equation}
\begin{equation}
\omega _{\tilde{f}_{m}}(h)\leq 2\omega (h)+\frac{h}{t},  \label{n12}
\end{equation}
and
\begin{equation}
\int \tilde{f}_{m}d\Psi _{n_{m}}\geq \delta -5\omega (t).  \label{n2}
\end{equation}

\item[(e)]  The non-decreasing sequence $\{g_{1}(x)\leq g_{2}(x)\leq ...\}$
such that $g_{m}(x)=g_{m-1}(x)+\tilde{f}_{m,L_{m}}^{(\varepsilon /3)}(x)$,
\begin{equation}
\left| \int g_{m-1}(x)d\Psi _{n_{m}}\right| \leq t,  \label{n22}
\end{equation}
$g_{m}(x)=0$ for $x\not{\in}K_{m}^{\varepsilon /3}$, and
\[
\omega _{g_{_{m}}}(h)\leq 2\omega (h)+\frac{h}{t}+\frac{3h}{\varepsilon }.
\]
\end{description}

Let
\[
g(x)=\lim_{m\rightarrow \infty }g_{m}(x).
\]
Clearly, $0\leq g(x)\leq 1$, and
\begin{equation}
\omega _{g}(h)\leq 2\omega (h)+\frac{h}{t}+\frac{3h}{\varepsilon }.
\label{n3}
\end{equation}

Since the numbers $\varepsilon $ and $t$ are fixed, the function $g\in
\mathcal{C}_\omega$. We show that, nevertheless, one can choose $\varepsilon
$ and $t$ such that $I_{n}:=\int gd\Psi _{n}\nrightarrow 0$.

To make exposition simpler, we replace the sequence $\{n_{m}\}$ by $\Bbb{N}$%
, write $\Psi _{n}$ instead of $\Psi _{n_{m}}$, and remove $\symbol{126}$
from $\tilde{f}$'s. All of this cannot cause misunderstanding.

By virtue of (\ref{n2}),
\[
I_{n}=\int f_{n}d\Psi _{n}+\int (g-f_{n})d\Psi _{n}\geq \delta -5\omega
(t)+\int (g-f_{n})d\Psi _{n}.
\]
For $J_{n}:=\int (g-f_{n})d\Psi _{n}$, we write
\[
|J_{n}|\leq \left| \int_{K_{n}}(g-f_{n})d\Psi _{n}\right| +|\Psi
_{n}|(K_{n}^{\complement })\leq \left| \int_{K_{n}}(g-f_{n})d\Psi
_{n}\right| +t.
\]
Now,
\[
\int_{K_{n}}(g-f_{n})d\Psi _{n}=\int_{L_{n}}+\int_{K_{n}\cap
K_{n-1}^{\varepsilon /3}}+\int_{K_{n}\setminus (L_{n}\cup
K_{n-1}^{\varepsilon /3})}:=J_{n1}+J_{n2}+J_{n3}.
\]
By construction, $J_{1n}=0$. For the second integral, we have
\begin{eqnarray}
|J_{n2}| &=&\left| \int_{K_{n-1}^{\varepsilon /3}}(g-f_{n})d\Psi
_{n}-\int_{K_{n}^{\complement }\cap K_{n-1}^{\varepsilon /3}}(g-f_{n})d\Psi
_{n}\right|  \nonumber \\
&\leq &\left| \int_{K_{n-1}^{\varepsilon /3}}(g-f_{n})d\Psi _{n}\right|
+\int_{K_{n}^{\complement }}|\Psi _{n}|(dx)  \nonumber \\
&\leq &\left| \int_{K_{n-1}^{\varepsilon /3}}gd\Psi _{n}\right| +\left|
\int_{K_{n-1}^{\varepsilon /3}}f_{n}d\Psi _{n}\right| +t,  \label{n4}
\end{eqnarray}
in view of (\ref{n0}).

By construction, $g(x)=g_{n-1}(x)$ for $x\in K_{n-1}^{\varepsilon /3}$. So,
\begin{equation}
\left| \int_{K_{n-1}^{\varepsilon /3}}gd\Psi _{n}\right| =\left|
\int_{K_{n-1}^{\varepsilon /3}}g_{n-1}d\Psi _{n}\right| =\left| \int
g_{n-1}d\Psi _{n}\right| \leq t,  \label{n5}
\end{equation}
by virtue of (\ref{n22}).

In view of (\ref{n11}) and (\ref{n12}),
\begin{eqnarray*}
\left| \int_{K_{n-1}^{\varepsilon /3}}f_{n}d\Psi _{n}\right| &=&\left|
\int_{K_{n-1}^{\varepsilon /3}\setminus K_{n-1}}f_{n}d\Psi _{n}\right| \leq
\omega _{f_{n}}(\varepsilon /3)\int \left| \Psi _{n}\right| (dx) \\
&\leq &\left( 2\omega (\varepsilon /3)+\frac{\varepsilon /3}{t}\right)
2=4\omega (\varepsilon /3)+\frac{2\varepsilon }{3t}.
\end{eqnarray*}
Thus,
\begin{equation}
|J_{n2}|\leq 4\omega (\varepsilon /3)+\frac{2\varepsilon }{3t}+2t.
\label{n6}
\end{equation}
To evaluate $J_{n3}$, first note that if $x\in K_{n}\setminus (L_{n}\cup
K_{n-1}^{\varepsilon /3})$, then $x\in K_{n-1}^{\varepsilon }$, and hence
\begin{eqnarray}
\left| \int_{K_{n}\setminus (L_{n}\cup K_{n-1}^{\varepsilon
/3})}f_{n}(x)d\Psi _{n}\right| &\leq &\left| \int_{K_{n}\setminus (L_{n}\cup
K_{n-1}^{\varepsilon /3})}\omega _{f_{n}}(\varepsilon )|\Psi _{n}|(dx)\right|
\nonumber \\
&\leq &\left( 2\omega (\varepsilon )+\frac{\varepsilon }{t}\right) 2=4\omega
(\varepsilon )+\frac{2\varepsilon }{t}.  \label{n7}
\end{eqnarray}
Now, let us observe that if $x\in K_{n}\setminus (L_{n}\cup
K_{n-1}^{\varepsilon /3})$ and $x\notin L_{n}^{\varepsilon /3}$, then $%
g(x)=0 $. So,
\[
\int_{K_{n}\setminus (L_{n}\cup K_{n-1}^{\varepsilon /3})}g(x)d\Psi
_{n}=\int_{(K_{n}\setminus L_{n})\cap L_{n}^{\varepsilon /3}}g(x)d\Psi _{n}.
\]
On the other hand, $(K_{n}\setminus L_{n})\cap L_{n}^{\varepsilon
/3}\subseteq K_{n-1}^{\varepsilon }$, and $g(x)\leq f_{n}(x)$ on $%
(K_{n}\setminus L_{n})\cap L_{n}^{\varepsilon /3}$. Thus, for the function $%
g $, we have the bound similar to (\ref{n7}), and
\[
|J_{n3}|\leq 8\omega (\varepsilon )+\frac{4\varepsilon }{t}.
\]

Combining the bounds above, we have
\[
|J_{n}|\leq 4\omega (\varepsilon /3)+8\omega (\varepsilon )+\frac{%
14\varepsilon }{3t}+3t\leq 12\omega (\varepsilon )+5\frac{\varepsilon }{t}%
+3t,
\]
and
\[
|I_{n}|\geq \delta -5\omega (t)-12\omega (\varepsilon )-5\frac{\varepsilon }{%
t}-3t\geq \delta -17\omega (t)-5\frac{\varepsilon }{t}-3t,
\]
because we choose $\varepsilon \leq t$. Without loss of generality we can
take $t<1$. Let $\varepsilon =t^{2}$. Then
\[
|I_{n}|\geq \delta -17\omega (t)-8t.
\]
Clearly, one can choose $t$ for which $|I_{n}|\geq \frac{\delta }{2}$ for
all $n$. $\blacksquare $

\end{document}